\documentstyle[11pt]{article}

\begin{document}

\title{Countable support iterations and large continuum}
\author{
Chaz Schlindwein \\
Division of Mathematics and Computer Science \\
Lander University \\
Greenwood, South Carolina 29649, USA\\
{\tt chaz@lander.edu}}

\maketitle

\def\implies{\Rightarrow}
\def\rk{{\rm rk}}
\def\overb{{\overline b}}
\def\overw{{\overline w}}
\def\overx{{\overline x}}
\def\overy{{\overline y}}
\def\overz{{\overline z}}
\def\hht{{\rm ht}}
\def\forces{\mathbin{\parallel\mkern-9mu-}}
\def\notforces{\,\nobreak\not\nobreak\!\nobreak\forces}
\def\Gen{{\rm Gen}}
\def\calm{{\cal M}}
\def\bfone{{\bf 1}}

\def\restr{\,\hbox{\vrule height8pt width.4pt depth0pt
\vrule height7.75pt width0.3pt depth-7.5pt\hskip-.2pt
\vrule height7.5pt width0.3pt depth-7.25pt\hskip-.2pt
\vrule height7.25pt width0.3pt depth-7pt\hskip-.2pt
\vrule height7pt width0.3pt depth-6.75pt\hskip-.2pt
\vrule height6.75pt width0.3pt depth-6.5pt\hskip-.2pt
\vrule height6.5pt width0.3pt depth-6.25pt\hskip-.2pt
\vrule height6.25pt width0.3pt depth-6pt\hskip-.2pt
\vrule height6pt width0.3pt depth-5.75pt\hskip-.2pt
\vrule height5.75pt width0.3pt depth-5.5pt\hskip-.2pt
\vrule height5.5pt width0.3pt depth-5.25pt}\,}

\def\overtau{{\overline\tau}}

\def\overt{{\overline t}}
\def\cali{{\cal I}}
\def\cf{{\rm cf}}
\def\sup{{\rm sup}}
\def\supt{{\rm supt}}
\def\dom{{\rm dom}}
\def\range{{\rm range}}
\def\calj{{\cal J}}
\def\calc{{\cal C}}

\centerline{{\bf Abstract}}

We prove that any countable support iteration formed with
posets with $\omega_2$-p.i.c.\  has $\omega_2$-c.c., assuming
CH in the ground model. This improves earlier results of Shelah
by removing the restriction on the length of the iteration.
Thus, we solve the problem of obtaining a large
continuum via such forcing iterations.

\vfill\eject

Shelah [5, chapter VIII] introduces the notion of $\omega_2$-p.i.c.\
forcing. He shows that if $\langle P_\xi\,\colon\xi\leq\kappa\rangle$ is a
countable support forcing iteration based on
$\langle \dot Q_\xi\,\colon\xi<\kappa\rangle$ such that
$\dot Q_\xi$ has $\omega_2$-p.i.c.\  in $V[G_{P_\xi}]$ for
all $\xi<\kappa$, then if
CH holds in $V$ and $\kappa\leq\omega_2$ 
 then $P_\kappa$ has
$\omega_2$-c.c. Many familiar forcings satisfy
$\omega_2$-p.i.c., such as the forcings to add a Sacks real,
a Mathias real, a Laver real, and so forth. Also,
$\omega_2$-p.i.c.\  largely subsumes $\omega_2$-e.c.c.,
as demonstrated in [1, lemma 57], [2, lemma 24],
and [3, lemma 15]. In this paper we eliminate the
restriction on the length of the iteration.  Also, the
hypothesis of our theorem is that each
$\dot Q_\xi$ has weak $\omega_2$-p.i.c.\  and the iteration
does not collapse $\omega_1$, where weak
$\omega_2$-p.i.c.\  can be viewed as ``$\omega_2$-p.i.c.\
minus properness.''

The main idea which is required, aside from
[5, chapter VIII] of course, is the notion of
$\dot P^M_{\eta,\kappa}$ from [4]. The statement
$p\forces_{P_\eta}``\dot q\in\dot P_{\eta,\kappa}\cap M[G_{P_\eta}]
$'' means that for every $p_1\leq p$ there is
$p_2\leq p_1$ and $q_1\in P_\kappa$ and $x\in M$ such that
$x$ is a $P_\eta$-name and $p_2\forces``\dot q=x=q_1\restr[\eta,
\kappa)$.'' In contrast, the statement $p\forces_{P_\eta}``\dot q\in
\dot P^M_{\eta,\kappa}$'' means that for every $p_1\leq p$
there is $p_2\leq p_1$ and $q_1\in P_\kappa\cap M$
such that $p_2\forces``\dot q= q_1\restr[\eta,\kappa)$.''
As always, the notation ``$q_1\restr[\eta,\kappa)$'' does not
refer to the check (with respect to $P_\eta$) of the restriction
of $q_1$ to the interval $[\eta,\kappa)$, but rather to the
$P_\eta$-name which is forced to be a function with
domain $[\check\eta,\check\kappa)$ such that for every
$\gamma$ in the interval $[\eta,\kappa)$ we have
that $q_1\restr[\eta,\kappa)(\check\gamma)$ is the $P_\eta$-name
for the $\dot P_{\eta,\gamma}$-name corresponding to the
$P_\gamma$-name $q_1(\gamma)$ (see [1, section 3] for
greater detail on this point). The notion of
$p\forces``\dot q\in\dot P^M_{\eta,\kappa}$'' is
exploited in [4] to prove preservation of semi-properness under
countable support (CS) iteration, preservation
of hemi-properness under CS iteration, and a theorem
giving a weak but sufficient condition for a CS
iteration to add no reals.  The main property of
$\dot P^M_{\eta,\kappa}$ which is
needed in [4] and in the present paper is the fact that
$\bfone\forces_{P_\eta}``(\forall\dot q\in\dot P_{\eta,\kappa}^M)
(\supt(\dot q)\subseteq\check M)$.'' The fact that
this holds is clear from the characterization of
$p\forces_{P_\eta}``\dot q\in\dot P^M_{\eta,\kappa}$''
given above; in any case a detailed proof is given
in [4, lemma 3].

\proclaim Definition 1. 
We say that\/ $(P,M,N,i,j)$ is\/ {\rm embryonic}
iff for some  sufficiently
large regular\/ $\lambda$ we have that\/ $M$ and\/ $N$ are
 countable elementary substructures of\/ $H_\lambda$ and\/
 $\omega_1<i<j<\omega_2$ and\/
 $\cf(i)=\cf(j)=\omega_1$ and\/
 $P\in M\cap N$ and\/ $i\in M$ and\/ $j\in N$ and\/
 $\sup(\omega_2\cap M)<j$ and\/
 $i\cap M=j\cap N$.

\proclaim Definition 2. We say that $(P,M,N,i,j,h)$ is {\rm
passable} iff
for some sufficiently large regular cardinal $\lambda$ we have
that $(P,M,N,i,j)$ is embryonic and
$h$ is an isomorphism from $M$ onto $N$ and
$h$ is the identity on $M\cap N$ and $h(i)=j$.

\proclaim Definition 3. We say that $P$ has {\rm weak $\omega_2$-p.i.c.}
iff
whenever $(P,M,N,\allowbreak
i,j,h)$ is passable and $p\in P\cap M$ then
there is
some $q\leq p$ such that $q\leq h(p)$.

\proclaim Definition 4.
Suppose $\langle P_\xi\,\colon\xi\leq\kappa\rangle$ is a countable
support iteration. We say that\/
$(P_\kappa,M,N,i,j,\eta,h,p)$ is\/ {\rm strictly
passable} iff\/ $(P_\kappa,M,N,i,j)$ is embryonic
and\/ $\eta\in\kappa\cap M\cap N$ and
$p\in P_\eta$ and $p\forces``h$ {\rm is an isomorphism
from $M[G_{P_\eta}]$ onto $N[G_{P_\eta}]$ and
the restriction of $h$ to $M[G_{P_\eta}]\cap N[G_{P_\eta}]$
is the identity and $h(i)=j$ and $\check N$ is the image
of $\check M$ under $h$.''}

\proclaim Definition 5. We say $P_\kappa$ is\/ {\rm
strictly weak $\omega_2$-p.i.c.} iff
whenever $(P_\kappa,M,N,\allowbreak i,
\allowbreak j,\allowbreak \eta,h,p)$ is strictly passable
and $p\forces_{P_\eta}``\dot q\in\dot P^M_{\eta,\kappa}$''
then there is $r\in P_\kappa$ such that
$r\restr\eta=p$ and $p\forces``r\restr[\eta,\kappa)\leq\dot q$
{\rm and $r\restr[\eta,\kappa)\leq h(\dot q)$.''}

\proclaim Lemma 6. Suppose $\langle P_\xi\,\colon
\xi\leq\kappa\rangle$ is a countable support iteration
based on $\langle \dot Q_\xi\,\colon\xi<\kappa\rangle$
and $P_\xi$ is strictly weak\/ $\omega_2$-p.i.c.\   whenever
$\xi<\kappa$. Suppose that if $\kappa=\gamma+1$ then $\bfone\forces_{
P_\gamma}``
\dot Q_\gamma$ {\rm has weak $\omega_2$-p.i.c.\  and
$\omega_1^{V_{\null_{\null}}}=\omega_1^{V[G_{P_\gamma}]}$.''}
 Then $P_\kappa$ is strictly weak\/ $\omega_2$-p.i.c.

Proof: Suppose $(P_\kappa,M,N,i,j,\eta,h,p)$ is strictly passable
and $p\forces``\dot q\in\dot P_{\eta,\kappa}^M$.''

Case 1: $\kappa=\gamma+1$.

We have that $(P_\gamma,M,N,i,j,\eta,h,p)$ is strictly passable,
so we may take $r_0\in P_\gamma$ such that
$r_0\restr\eta=p$ and $p\forces``r_0\restr[\eta,\gamma)\leq
\dot q\restr\gamma$ and $r_0\restr[\eta,\gamma)\leq
h(\dot q\restr\gamma)$.''
Now, the restriction of $h$ to $M[G_{P_\eta}]^{\dot P_{\eta,\gamma}}$,
which by notational convention is the set of all
$\dot P_{\eta,\gamma}$-names in $M[G_{P_\eta}]$,
induces an isomorphism $h^*\in V[G_{P_\gamma}]$ from
$M[G_{P_\gamma}]$ onto $N[G_{P_\gamma}]$ such that
$h^*(i)=j$ and the restriction of $h^*$ to
$M[G_{P_\gamma}]\cap N[G_{P_\gamma}]$ is the identity and $\check N$
is the image of $\check M$.
Hence we may use the fact that
$\dot Q_\gamma$ has weak $\omega_2$-p.i.c.\  to
take $\dot r_1\in\dot Q_\gamma$
such that $r_0\forces
``\dot r_1\leq\dot q(\gamma)$ and $\dot r_1\leq h^*(\dot q(\gamma))$.''
Then let $r=(r_0,\dot r_1)\in P_\kappa$.
We have that $r$ is as required.

Case 2: $\kappa$ is a limit ordinal.

Let $\alpha=\sup(\kappa\cap M\cap N)$ and take
$\langle \alpha_n\,\colon n\in\omega\rangle$ an increasing
sequence of ordinals from $\alpha\cap M\cap N$ cofinal
in $\alpha$ such that $\alpha_0=\eta$.
Build $\langle p_n\,\colon n\in\omega\rangle$ such that
$p_0=p$ and for every $n\in\omega$ we have
$p_{n+1}\restr\alpha_n=p_n$ and $p_n\forces``p_{n+1}\restr
[\alpha_n,\alpha_{n+1})\leq\dot q\restr[\alpha_n,\alpha_{n+1})$
and $p_{n+1}\restr[\alpha_n,\alpha_{n+1})\leq
h^*_n(\dot q\restr[\alpha_n,\alpha_{n+1}))$''
where $h^*_n$ is the isomorphism from
$M[G_{P_{\alpha_n}}]$ onto $N[G_{P_{\alpha_n}}]$ induced
by $h$ as in the successor case above.
This is possible by the induction hypothesis and
[4, lemma 5].
Take $r\in P_\kappa$ such that $r\restr\alpha_n=p_n$
for every $n\in\omega$, and $r(\xi)=\dot q(\xi)$ for
all $\xi\in\kappa\cap M$ such that $\alpha\leq\xi$,
and $r(\xi)=h(\dot q)(\xi)$ for all $\xi\in\kappa\cap N$
such that $\alpha\leq\xi$, and
$r(\xi)=\bfone_{\dot Q_\xi}$ in all other cases.
Then $r$ is as required.
As noted earlier, the main point is that $p\forces``\supt(\dot q)
\subseteq\check M$.'' 

We repeat [5, VIII.2.3] (see also [1, lemma 42] but note
by the way that the proof of [1, lemma 43] is incorrect;
the present paper is in part a repair of this deficiency).

\proclaim  Lemma 7.  Suppose CH holds and $P$ has
weak $\omega_2$-p.i.c. Then $P$ has $\omega_2$-c.c.

Proof: Take $\lambda$ a sufficiently
large regular cardinal.
Given $\langle p_i\,\colon\cf(i)=\omega_1$ and
$i<\omega_2\rangle$ (potentially, a counterexample
to $\omega_2$-p.i.c.), take for each such $i$ a
countable elementary submodel $N_i$ of $H_\lambda$
such that $\{p_i,i,P\}\subseteq N_i$.
It suffices to show that there are $\eta<\xi<\omega_2$ and $h$
such that $\eta\cap N_\eta=\xi\cap N_\xi$ and $\sup(\omega_2\cap
N_\eta)<\xi$ and $h$ is an isomorphism from $N_\eta$ onto $N_\xi$
and $h(\eta)=\xi$ and the restriction of $h$ to $N_\eta\cap
N_\xi$ is the identity. Take $f(i)=\sup(i\cap N_i)$
for all $i<\omega_2$ such that $\cf(i)=\omega_1$.
Take $S_0\subseteq\{i<\omega_2\,\colon\cf(i)=
\omega_1\}$ stationary and $\gamma<\omega_2$
such that $(\forall i\in S_0)(f(i)=\gamma)$. By CH we may take
 $S_1\subseteq S_0$ such that $\vert S_1\vert = \aleph_2$
 and whenever
$i<j$ are both in $S_1$ then $i\cap N_i=j\cap N_j$.
Take $S_2\subseteq S_1$ of size $\aleph_2$ such that
whenever $i<j$ are both in $S_2$ then
$\sup(\omega_2\cap N_i)<j$. By $\Delta$-system, take
$S_3\subseteq S_2$ of size $\aleph_2$ and a countable $N$ such
that whenever $i$ and $j$ are distinct elements of $S_3$
then $N_i\cap N_j=N$. Fix an enumeration
$\langle c_k\,\colon k\in\omega\rangle$ of $N$.
Let $N_i^+=\langle N_i;\in,i,c_0,c_1,\ldots\rangle$.
Up to isomorphism there are only $\aleph_1$-many possible
$N^+_i$, so the lemma is established.

Thus we have proved the following:

\proclaim Theorem 8. Suppose $\langle P_\xi\,\colon\xi\leq
\kappa\rangle$ is a countable support iteration based
on\/ $\langle\dot Q_\xi\,\colon\xi<\kappa\rangle$
and for each $\xi<\kappa$ we have that
$\dot Q_\xi$ has weak $\omega_2$-p.i.c.\  in $V[G_{P_\xi}]$.
Suppose also that $\omega_1^{V_{\null_{\null}}}=\omega_1^{V[G_{P_\kappa}]}$ and
that\/ {\rm CH} holds in $V$. Then $P_\kappa$ has
$\omega_2$-c.c.

We also have the following:

\proclaim Fact 9. Suppose $\langle P_\xi\,\colon\xi\leq\kappa\rangle$
is a countable support forcing iteration. Then no reals are added at limit
stages of uncountable cofinality.

Proof: Suppose $\alpha\leq\kappa$ and $\cf(\alpha)>\omega$ and
$q\forces_{P_\alpha}``\dot r\in{}^\omega\omega$.'' Take $M$
a countable model containing all relevant data. Let $\beta=
\sup(\alpha\cap M)$ and let $\langle\beta_n\,\colon n\in\omega\rangle$
be an increasing sequence of ordinals from $\beta\cap M$ cofinal
in $\beta$ with $\beta_0=0$. Build
$\langle p_n,q_n\,\colon n\in\omega\rangle$ such that $q=q_0$ and
$p_n\forces_{P_{\beta_n}}``q_n\in \dot P_{\beta_n,\alpha}^M$ and
$q_n\leq q_{n-1}\restr[\beta_n,\alpha)$ and $q_n$ decides
the value of $r(n-1)$'' and $p_{n+1}\restr\beta_n=p_n$ and
$p_n\forces``p_{n+1}\restr[\beta_n,\beta_{n+1})\leq
q_n\restr\beta_{n+1}$.'' Then take $q'\in P_\alpha$
such that $\supt(q')\subseteq\beta$ and $q'\restr\beta_n=p_n$
for all $n\in\omega$. We have that $q'\leq q$ and
$q'\restr\beta\forces``\bfone\forces_{\dot P_{\beta,\alpha}}`
r=\check s$'\thinspace'' for some $s\in V[G_{P_\beta}]$.

\proclaim Fact 10.  The ${}^\omega\omega$-bounding property,
the Sacks property, the Laver property, etc., are preserved by
countable support iteration (without the assumption of properness).

This follows from Fact 9, together with the arguments of [5, section VI.2];
the only place [5] uses the properness assumption is to handle
the uncountable cofinality case.

\vfill\eject

\noindent{\bf References}

\medskip

[1] Schlindwein, C.,
``Consistency of Suslin's hypothesis, a non-special Aronszajn tree,
and GCH,'' {\sl Journal of Symbolic Logic,} vol.~{\bf 59}, pp.~1--29, 1994

                  \medskip

[2] Schlindwein, C., ``Suslin's hypothesis does not imply
stationary antichains,'' {\sl Annals of Pure and Applied
Logic,} vol.~{\bf 64}, pp.~153--167, 1993

\medskip

[3] Schlindwein, C., ``Countably paracompact Aronszajn trees
and CH.''

\medskip

[4] Schlindwein, C., ``Preservation theorems for
countable support forcing iterations.''

                          \medskip

[5] Shelah, S., {\sl Proper Forcing,} Lecture Notes in Mathematics
{\bf 940}, Springer-Verlag, 1982

\vfill\eject

\end{document}